\documentclass[12pt]{amsart}
\usepackage{amsmath,latexsym,amssymb}
\usepackage{amscd}

\input xypic
\newcommand{\rar}{\rightarrow}
\newcommand{\lar}{\longrightarrow}

\newtheorem{theorem}{Theorem}[section]
\newtheorem{lemma}[theorem]{Lemma}
\newtheorem{corollary}[theorem]{Corollary}
\newtheorem{proposition}[theorem]{Proposition}

\newtheorem{remark}[theorem]{Remark}
\newtheorem{example}[theorem]{Example}

\newtheorem{problem}[theorem]{Problem}

\def\Ass{\mbox{\rm Ass}}
\def\demo{\noindent{\bf Proof. }}
\def\depth{\mbox{\rm depth}}
\def\ds{\displaystyle}
\def\Ext{\mbox{\rm Ext}}
\def\gr{\mbox{\rm gr}}
\def\height{\mbox{\rm height}}
\def\Hom{\mbox{\rm Hom}}
\def\l{\lambda}
\def\m{\mathfrak{m}}

\def\n{\mathfrak{n}}
\def\p{\mathfrak{p}}
\def\q{\mathfrak{q}}
\def\a{\mathfrak{a}}
\def\QED{\hfill$\Box$}
\def\Rees{\mathcal{R}}
\def\supp{\mbox{\rm supp }}
\def\Tor{\mbox{\rm Tor}}

\def\XX{{\bf X}}
\def\xx{{\bf x}}

\def\BB{{\bf B}}

\begin{document}

\title[The Signature of the  Chern Coefficients of Local Rings
]{The Signature of the Chern Coefficients \\ of Local Rings}

\thanks{AMS 2000 {\em Mathematics Subject Classification}.
Primary 13A30; Secondary 13B22, 13H10, 13H15.}

\author{Laura Ghezzi}
\address{Department of Mathematics, New York City College of
Technology-Cuny}
\address{300 Jay Street, Brooklyn, NY 11201,  U.S.A.}
\email{lghezzi@citytech.cuny.edu}

\author{Jooyoun Hong}
\address{Department of Mathematics, Southern Connecticut State
University}
\address{501 Crescent Street, New Haven, CT 06515-1533,  U.S.A.}
\email{hongj2@southernct.edu}

\author{Wolmer V. Vasconcelos}
\thanks{The first author is partially supported by a grant from the City University of New York PSC-CUNY Research Award Program-39. The third author is partially supported by the NSF}
\address{Department of Mathematics, Rutgers University}
\address{110 Frelinghuysen Rd, Piscataway, NJ 08854-8019, U.S.A.}
\email{vasconce@math.rutgers.edu}

\begin{abstract}
This paper considers
 the following conjecture: If $R$ is
 an unmixed, equidimensional local ring that is a
homomorphic image of a Cohen-Macaulay local ring, then for any
ideal $J$  generated by a system of parameters, the Chern
coefficient
  $e_1(J)< 0$ is equivalent to $R$ being non Cohen-Macaulay.
The conjecture is established if $R$ is a homomorphic image of a
Gorenstein ring, and for all universally catenary integral domains
containing fields. Criteria for
the detection of Cohen-Macaulayness in equi-generated graded
modules are derived.

\end{abstract}

\maketitle

\section{Introduction}

Let $(R, \mathfrak{m})$  be a Noetherian local ring of dimension $d>0$, and let $I$ be an
$\mathfrak{m}$-primary ideal. We will consider the
set of multiplicative, decreasing  filtrations of $R$ ideals, $\BB=\{ I_n,
 I_0=R,  I_{n+1}=I I_n, n\gg 0 \}$, integral over the $I$-adic filtration. They are
conveniently coded in the corresponding Rees algebra and its associated graded ring
\[ \Rees(\BB) = \sum_{n\geq 0} I_nt^n, \quad
\gr_{\BB}(R) = \sum_{n\geq 0} I_n/I_{n+1}.
 \] One of our goals is to study
cohomological properties of
these filtrations.
For that we will focus on the role of the Hilbert polynomial of the Hilbert--Samuel
function $\lambda(R/I_{n+1})$,
\[ H_{\BB}^1(n) =P_{\BB}^1(n) \equiv \sum_{i=0}^{d} (-1)^i
e_i(\BB)
{{n+d-i}\choose{d-i}}, \] particularly of the multiplicity,
$e_0(\BB)$, and the {\em Chern  coefficient},
 $e_1(\BB)$. For Cohen-Macaulay rings,
many penetrating relationships between these coefficients
 have been proved,
beginning with Northcott's \cite{Northcott}. More recently, similar
questions have been examined in general Noetherian local rings and
among those pertinent to our concerns are \cite{GoNi03}, \cite{RV08}
and \cite{chern}.

Here we extend several of the results of \cite{chern} on
 the meaning of the sign of $e_1(\BB)$,
particularly in the case of $I$-adic filtrations. Our main
 results are centered around the following question:

\medskip

Let $J$ be an ideal generated by a system of parameters. Under
which conditions is $e_1(J)< 0$ equivalent to $R$ being not
Cohen-Macaulay?  We conjecture that this is so
whenever $R$ is an unmixed, equidimensional local ring that is a
homomorphic image of a Cohen-Macaulay local ring.

\medskip

There are reasons for the interest in these numbers. To make the
discussion more direct, we assume that the residue field of $R$ is
infinite.

\medskip

(1)
 Clarifying the role of the sign of $e_1(J)$ in the
Cohen-Macaulay property of $R$ could be used as a {\em scale} to
classify non-CM rings.

\medskip

 (2) In the study of the normalization $\BB$ of $R[Jt]$, the
expression (see \cite{ni1}, \cite{chern})
\[ e_1(\BB)-e_1(J)\]
bounds the length of certain computations in the construction of
normalizations.  It highlights the need to look for upper bounds of
$e_1(\BB)$ and for lower bounds of $e_1(J)$.

\medskip

(3) The value of $e_1(J)$ occurs as a correction term in the
extensions of several well-known formulas in the theory of the
Hilbert polynomials. We highlight two of them. A classical
result of Northcott (\cite{Northcott}) asserts that if $R$ is
Cohen-Macaulay, then
\[ e_1(I)\geq e_0(I)-\lambda(R/I).\]
For arbitrary Noetherian rings, if $J$ is a minimal reduction of $I$,
 Goto and
Nishida (\cite[Theorem 3.1]{GoNi03})proved that
\[ e_1(I)-e_1(J)\geq e_0(I)-\lambda(R/I),\] which gives
\[ e_1(\BB)-e_1(J)\geq e_0(I)-\lambda(R/\overline{I}),\] since
$e_1(\BB)\geq e_1(I)$ for all such ideals $I$. It is a
 formula which is relevant to a conjecture of \cite{chern}, on whether
$e_1(\BB)\geq 0$.

A different kind of relationship  given by Huckaba and Marley
(\cite[Theorem 4.7]{HuM}) for Cohen-Macaulay rings,
\[e_1(\BB) \leq \sum_{n\geq 1}\lambda(I_n/JI_{n-1}),\]
is extended by Rossi and Valla (\cite[Theorem 2.11]{RV08}) to general
filtrations to an expression that replaces $e_1(\BB)$ in the
inequality above by $e_1(\BB)-e_1(J)$.

\medskip

In our main result (Theorem~\ref{e1dMCM})
we show that the above Conjecture holds if $R$ is a homomorphic image of a Gorenstein
ring, or milder extensions that allow an embedding of $R$ into a
(small) Cohen-Macaulay module over a possibly larger ring.
The proof is a variation of an argument
in \cite{chern}, but turned more abstract. Another result
establishes the Conjecture for universally catenary local
domains containing a field (Theorem~\ref{e1Hochster}).

The same question can be asked about filtrations of modules regarding
the negativity of the coefficient $e_1$ of the corresponding
associated graded module. In case $R$ is a polynomial ring over a field
and $M$ is a graded, torsion--free $R$-module, this extension has a
surprising application to the Cohen-Macaulayness of $M$. In the
special case of modules generated in the same degree, the Hilbert
coefficient $e_1(M)$, which may be different from
$e_1(\gr_{\xx}(M))$, can alone decide whether $M$ is Cohen-Macaulay
or not (Corollary~\ref{e1vsCMa}) .

We thank Rodney Sharp and Santiago Zarzuela for sharing with us their
expertises on balanced big Cohen-Macaulay modules.
We are also grateful to Shiro Goto for sharing with us a sketch of his
proof of a more extended version of   Theorem~\ref{e1dMCM}. The
Yokohama Conference on Commutative Algebra, on March 2008, was the
setting where some of the questions arose.

\section{Preliminaries}

We will assemble quickly here some facts about associated graded
modules and their Hilbert functions. As a general reference for
unexplained terminology and basic results, we shall use \cite{BH}.

\medskip

Let $(R, \m)$ be a Noetherian local ring, $I$ an $\mathfrak{m}$-primary ideal, and $M$ a nonzero finite $R$--module of dimension $d$. The associated graded ring ${\ds \gr_I(R)=\bigoplus_{i=0}^{\infty} I^i/I^{i+1} }$ is a standard graded ring with $[\gr_I(R)]_0=R/I$ Artinian. The associated graded module ${\ds \gr_I(M)=\bigoplus_{i=0}^{\infty} I^iM/I^{i+1}M }$  of $I$ with respect to $M$ is a finitely generated graded $\gr_I(R)$--module. The {\em Hilbert--Samuel function} $\chi^I_M(n)$ of $M$ with respect to $I$ is
\[
\chi^I_M(n)=\l(M/I^{n+1}M) = \sum_{i=0}^n \l (I^iM/I^{i+1}M).
\] For sufficiently large $n$, the Hilbert--Samuel function $\chi^I_M(n)$ is of polynomial type :
\[
\chi^I_M(n)=\sum_{i=0}^{d} (-1)^i
e_i(I,M)
{{n+d-i}\choose{d-i}}. \]

For an $R$--module $M$ of finite length, we denote the length of $M$ by $\l(M)$.

\begin{lemma}\label{Prelim2}
Let $(R, \m)$ be a Noetherian local ring and let $I$ be an $\m$--primary ideal. Let $0 \rar T \rar M \rar N \rar 0$ be an exact sequence of finitely generated $R$--modules. Assume that $M$ has dimension $d\ge 2$ and that $T$ has finite length. Then $e_1(I, M)=e_1(I, N)$.
\end{lemma}

\demo From the following commutative diagram with exact rows

\[\begin{CD}
0 @>>> T \cap I^{n+1}M @>>> I^{n+1}M @>>> I^{n+1}N @>>> 0 \\
 & & @VVV  @VVV  @VVV \\
0 @>>> T  @>>> M @>>> N @>>> 0
\end{CD}\] we get an exact sequence
\[
0 \lar T/(T \cap I^{n+1}M) \lar M/I^{n+1}M \lar N/I^{n+1}N \lar 0.
\] By the Artin--Rees Lemma, $T \cap I^{n+1}M=0$ for all sufficiently large $n$. Hence
we get for all $n\gg 0$
\[
\l(T)-\l(M/I^{n+1}M)+\l(N/I^{n+1}N)=0.\]

Let $d'$ be the dimension of $N$.
There are Hilbert polynomials such that
\[
{\ds \l(M/I^{n+1}M)= \sum_{i=0}^d (-1)^i e_i(I,M){{d+n-i}\choose{d-i}}} \]
and
\[{\ds \l(N/I^{n+1}N)= \sum_{i=0}^{d'} (-1)^i e_i(I, N){{d'+n-i}\choose{d'-i}}
}\]for all sufficiently large $n$. Therefore $d=d'$,
\[
e_i(I,M)=e_i(I,N) \quad \mbox{\rm for all}\; i=0, \ldots, d-1 \quad \mbox{\rm and} \quad e_d(I,N)=e_d(I,M)+(-1)^{d+1}\l(T).
\] In particular, since $d \geq 2$, we have $e_1(I,M)=e_1(I,N)$. \QED

\bigskip

Let $R$ be a ring, $I$ an $R$--ideal and $M$ an $R$--module. An element $h \in I$ is called a {\em superficial element} of $I$ with respect to $M$ if there exists a positive integer $c$ such that
\[
(I^{n+1}M :_M h) \cap I^cM = I^n M
\] for all $n \geq c$.
A detailed discussion on superficial elements can be found in several sources, but we
especially benefited from the treatment in \cite[Theorem 1.5]{RV07}
of the construction of superficial elements. It depends simply on
showing that certain finitely generated modules cannot be written
as a union of a finite set of proper submodules. The existence of
such elements is  guaranteed if the residue field of
$R$ is infinite. Its usefulness for our purposes is expressed in the
following result.

\begin{proposition}\label{genhs} {\rm (\cite[(22.6)]{Nagata})} Let $(R, \m)$ be a Noetherian local ring, $I$ an $\mathfrak{m}$-primary ideal, and $M$ a nonzero finitely generated $R$--module of dimension $d$. Let $h$ be a superficial element of $I$ with respect to $M$. Then
the Hilbert coefficients of $M$ and $M/hM$ satisfy
\[
e_i(I,M) = \left\{
\begin{array}{ll}
e_i(I/(h), M/hM) &\quad \mbox{\rm for }\;\; i< d-1.\\
e_{d-1}(I/(h), M/hM) + (-1)^{d-1} \lambda(0:_M h). &\quad \mbox{\rm for }\;\; i=d-1.
\end{array}  \right. \]
\end{proposition}

\section{Cohen-Macaulayness versus the vanishing of the Euler number}

In this section we develop an  abstract approach to the
relationship between the signature of $e_1$ and the
Cohen-Macaulayness of a local ring (see \cite[Theorem 3.1]{chern}). We also give a more
general but still self-contained proof of the main result of
\cite{chern} that avoids the use of big
Cohen-Macaulay modules.

\begin{lemma}\label{lifting}
Let $(S, \n)$ be a Cohen--Macaulay local ring of dimension $d$ with infinite residue field and let $R=S/\p$, where $\p$ is a  minimal prime ideal of $S$.
Let $x_1, \ldots, x_d$ be a system of parameters of $R$. Then there exists a system of parameters $a_1, \ldots, a_d$ of $S$ such that $x_i=a_i+ \p$ for each $i$.
\end{lemma}

\demo Let $\m$ denote the maximal ideal of $R$ and let $\p=(c_1, \ldots, c_s)$. Let $x_1=b_1+\p$ for some $b_1 \in S$ and let $\p_1, \ldots, \p_n$ be the minimal primes of $S$ different from $\p$. We claim that there exists $\l \in S \setminus \n$ such that $b_1 + \l c_1 + \cdots + \l^s c_s \not \in \p_i$ for all $i=1, \ldots, n$. Suppose not: since $S/\n$ is infinite, there exist $\l_1, \ldots, \l_{s+1} \in S \setminus \n$ such that $\l_i + \n \neq \l_j +\n$ whenever $i \neq j$ and such that $b_1+ \l_ic_1+ \cdots+ \l_i^s c_s \in \p_k$ for some fixed $k$. Let $A$ be the Vandermonde matrix determined by $\l_i$, $1\leq i\leq s+1$. We have
\[
A \left[ \begin{array}{c} b_1 \\ c_1 \\ \vdots \\ c_s  \end{array}\right]=\left[
\begin{array}{ccccc}
1 & \l_1 & \l_1^2 & \cdots & \l_1^s \\
1 & \l_2 & \l_2^2 & \cdots & \l_2^s \\
\vdots & \vdots &\vdots &\ddots &\vdots \\
1 & \l_{s+1} & \l_{s+1}^2 & \cdots & \l_{s+1}^s \\
\end{array}
\right] \left[ \begin{array}{c} b_1 \\ c_1 \\ \vdots \\ c_s  \end{array}\right]=
\left[\begin{array}{c} g_1 \\ g_2 \\ \vdots \\ g_{s+1} \end{array}\right],
\]where $g_1, \ldots, g_{s+1} \in \p_k$. Since $A$ is invertible,
$c_1, \ldots, c_s \in \p_k$ so that $\p \subseteq \p_k$, a
contradiction. Let $a_1=b_1 + \l c_1 + \cdots + \l^s c_s$, with $\l
\in S \setminus \n$, be such that $a_1 \not \in \p_i$ for all $i=1,
\ldots, n$.  Hence $a_1$ is not contained in any minimal prime ideal of $S$ and $a_1+\p=x_1$.

\medskip

Let $\q_1, \ldots, \q_m$ be the minimal primes of $a_1S$ that do not
contain $\p$. Let $x_2=b_2+\p$ for some $b_2 \in S$. Then similarly
as shown above, there exists $\tau \in S \setminus \n$ such that
$b_2+\tau c_1 + \cdots + \tau^s c_s \not \in \q_i$ for all $i$. Let $a_2=b_2+\tau c_1 + \cdots + \tau^s c_s$. Then $a_2$ is not contained in any minimal prime ideal of $a_1S$ and $a_2+\p=x_2$.
 Now inductively we show that there exists
 a system of parameters $a_1, \ldots, a_d$ of $S$ such that $a_i+\p=x_i$ for all $i$. \QED

\begin{lemma}\label{H1}
Let $(S, \n)$ be a Cohen--Macaulay complete local ring of dimension $d \geq 2$ and let $M$ be a finitely generated $S$--module of dimension $d$ satisfying Serre's condition {\rm ($S_1$)}.
 Then $H_{\n}^1(M)$ is a finitely generated $S$--module.
\end{lemma}

\demo Let $k$ denote the residue field of $S$ and $\omega$ the canonical module of $S$. Since $M$ satisfies Serre's condition {\rm ($S_1$)}, $\Ext^{d-1}_S(M, \omega)$ is zero at each localization at $\p$ such that $\p \neq \n$. Thus it is a module of finite length. By Grothendieck duality (\cite[3.5.8]{BH}), we have
\[H_{\n}^1(M) \simeq \Hom_S(\Ext_S^{d-1}(M, \omega), E(k)).\] By Matlis duality(\cite[3.2.13]{BH}), $H_{\n}^1(M)$ is finitely generated. \QED

\bigskip

Let $(R, \mathfrak{m})$ be a Noetherian local ring of dimension
$d\geq 2$. The enabling idea is the embedding of $R$ into a
Cohen-Macaulay (possibly big Cohen-Macaulay) module $E$,
\[ 0 \rar R \lar E \lar C \rar 0.\]
Unfortunately, it may not be always possible to find the appropriate
$R$-module $E$. Instead we will seek embed $R$ into a Cohen-Macaulay
module $E$ over a ring $S$ closely related to $R$ for the purpose of
computing associated graded rings of adic-filtrations.
The following is our main result.

\begin{theorem} \label{e1dMCM}
Let $(R, \m)$ be a Noetherian local domain of dimension
$d\geq 2$, which is a homomorphic image of a Cohen-Macaulay local ring $(S, \n)$. If $R$ is not Cohen-Macaulay, then $e_1(J)<
0$ for any $R$--ideal $J$ generated by a system of parameters.
\end{theorem}

\demo Let $R=S/\p$. We may assume that $S$ has
 infinite residue field. If $\height(\p) \geq 1$, we replace $S$ by
$S/L$, where $L$ is the $S$--ideal generated by a maximal regular sequence in $\p$.
This means that we may assume that $\dim R=\dim S$, and that
$\mathfrak{p}$ is a minimal prime of $S$. In particular, we have an
exact sequence of $S$-modules
\[ 0 \lar R \lar S \lar C \lar 0.\]
Let $J=(x_1, \ldots, x_d)$ be an $R$--ideal generated by a system of parameters. Then by Lemma~\ref{lifting}, there exists a system of parameters $a_1, \ldots, a_d$ of $S$ such that $x_i=a_i+\p$ for each $i$. Let $I=(a_1, \ldots, a_d)S$. Since $IR=J$, the associated graded ring $\gr_I(R)$ of $I$ with respect to the $S$--module $R$ is equal to the associated graded ring $\gr_J(R)$ of $J$. In particular, $e_1(I,R)=e_1(J)$.
Therefore, for the purpose of constructing  $e_1(J)$, we treat $R$ as
an $S$-{\em module} and use the $I$-adic filtration.

\medskip

Now proceed as in the proof of \cite[Theorem~3.1]{chern}. Let $a$ be a superficial element for $I$ with respect to the $S$--module $R$  which is not contained in any associated prime of $C$ distinct from
$\n$. We may assume that $a=a_1$. Reduction modulo $aS$ gives rise to the exact
sequence
\[ 0 \rar T = \Tor_1^S(S/aS, C) \lar  R/aR \lar S/aS \lar C/aC \rar
 0,\] where $T= (0:_C a)\subset C$. This shows that the associated
 primes of $T$ contain $a$. Therefore $T$ is either zero, or $T$ is a non-zero module of finite support.

\medskip

Let $R\,'=R/aR$, $I\,'=I/(a)$, and denote
the image of $R/aR$ in $S/aS $ by $S\,'$.

\medskip

Now we use induction on $d \geq 2$ to show that if $R$ is not Cohen--Macaulay, then $e_1(I,R) <0$.
Let $d=2$. Notice that $S\,'$ is a Cohen-Macaulay ring
of dimension $1$. We have that depth $C=0$, and so $T \neq 0$.
As in the proof of \cite[Theorem 3.1]{chern}, we obtain $e_1(I\,', R\,')=-\l(T)$.
Hence by Proposition~\ref{genhs}, $e_1(I,R)=-\l(T) <0$.

\medskip

Suppose $d >  2$.
Consider the exact sequence of $S/aS$--modules : $0 \rar T \rar R\,' \rar S\,' \rar 0$.
By Lemma~\ref{Prelim2}, we have $e_1(I', R\,')=e_1(I', S\,')$ since $\dim(S/aS)=d-1 \geq 2$.
Now it is enough to show that $S\,'$ is not a Cohen--Macaulay $S/aS$--module. Then since $\dim(S\,')=d-1$, by induction we get $e_1(I',S\,') <0$, and we conclude using Proposition~\ref{genhs}.

\medskip

\noindent Suppose that $S\,'$ is a Cohen--Macaulay $S/aS$--module. Let $\n$ denote the maximal ideal of $S/aS$ as well and let $H_{\n}^i( \cdot)$ denote the $i$th local cohomology. We are going to use the argument of \cite[Proposition 2.1]{HuL}. From the exact sequence $0 \rar T \rar R\,' \rar S\,' \rar 0$, we obtain a long exact sequence:
\[
0 \rar H_{\n}^0(T) \rar H_{\n}^0(R\,') \rar H_{\n}^0(S\,') \rar H_{\n}^1(T) \rar H_{\n}^1(R\,') \rar H_{\n}^1(S\,').
\]Since $S\,'$ is Cohen--Macaulay of dimension $d-1 \geq 2$, we have $H_{\n}^0(S\,')=0=H_{\n}^1(S\,')$. Since $T$ is a torsion module, we have $H_{\n}^0(T)=T$ and $H_{\n}^1(T)=0$. Therefore $T \simeq H_{\n}^0(R\,') $ and $H_{\n}^1(R\,')=0$. Now from the exact sequence of $S$--modules
\[0 \lar R \stackrel{\cdot a}{\lar} R \lar R/aR \lar 0 \]
we obtain the following exact sequence:
\[
0  \lar T \simeq H_{\n}^0(R\,') \lar H_{\n}^1(R) \stackrel{\cdot a}{\lar} H_{\n}^1(R) \lar H_{\n}^1(R\,')=0.
\] Therefore $H_{\n}^1(R)=aH_{\n}^1(R)$. Moreover, once $R$ is embedded in $S$, we may assume that $S$ is a complete
local ring. By Lemma~\ref{H1}, $H_{\n}^1(R)$ is finitely generated. By Nakayama Lemma, we have that $H_{\n}^1(R)=0$ so that $T=0$. It follows that $R/aR=R\,' \simeq S\,'$, where $S\,'$ is Cohen--Macaulay. This means that $R$ is Cohen--Macaulay, which is a contradiction. \QED

\medskip

\begin{remark}\label{e1remark} {\rm The proof of Theorem \ref{e1dMCM} can be extended from integral domains to more general local rings,
$R=S/L$, where $S$ is Cohen-Macaulay and $\dim R=\dim S$, if $R$ can be embedded into a maximal Cohen-Macaulay $S$-module. Notice that in order to embed $R$ into a maximal Cohen--Macaulay $S$-module at a minimum we need to require that $R$ be unmixed and equidimensional.}
\end{remark}

\medskip

There is room for the following problem:

\begin{problem} {\rm Let $R=S/L$ (unmixed and equidimensional as
above),   where $S$ is a Cohen-Macaulay local
ring
of dimension $\dim R$. Characterize those $R$ that can be embedded into
a maximal Cohen--Macaulay $S$-module. Note that $L$ may be assumed to be a primary
ideal: If
$L= \cap Q_i$ is a primary decomposition, we have the embedding
\[ S/L \hookrightarrow S/Q_1 \oplus \cdots \oplus S/Q_n.
 \]
}
\end{problem}

\medskip

We make  an elementary observation of what it takes to embed
  $R$ into a free $S$-module (see \cite[Theorem A.1]{Vas68}).

\begin{proposition}\label{gor} Let $S$ be a Noetherian  ring and $L$ an
ideal of codimension zero without embedded components. If $R=S/L$,
there is an embedding $R\rar F$ into a free
$S$--module $F$ if and only if $L=0:(0:L)$. In particular, this
condition always holds if the total ring of fractions of  $S$ is a Gorenstein ring.
\end{proposition}

\demo Let $\{a_1, \ldots, a_n\}$ be a generating set of  $0:L$, and
consider the   mapping $\varphi: S
\rar S^n$, $\varphi(1)=(a_1, \ldots, a_n)$; its kernel
is isomorphic to $0:(0:L)$. This shows that the equality $L=0:(0:L)$
is required for the asserted embedding.

Conversely, given an embedding $\varphi:S/L \rar S^n$, let $(a_1,
\ldots, a_n)\in S^n$ be the image of a generator of $S/L$. The
ideal $\a$ these entries generate  is annihilated by $L$, and so $\a\subset
0:L$. Since $ 0:\a=L$, we have $0:(0:L)\subset 0:\a= L$.

If the total ring of fractions of $S$ is Gorenstein, to prove that
 $0:(0:L)\subset  L$, it suffices to localize at the associated
 primes of $L$, all of which have codimension zero and a
 localization which is Gorenstein. But the double annihilator
 property is characteristic of such rings. \QED

\medskip

\begin{corollary} \label{e1dMCM2}
Let $(R, \m)$ be an unmixed and equidimensional Noetherian local ring of dimension
$d\geq 2$, which is a homomorphic image of a Gorenstein local ring. If $R$ is not Cohen-Macaulay, then $e_1(J)<
0$ for any $R$--ideal $J$ generated by a system of parameters.
\end{corollary}

\demo Let $R=S/L$. If $\height(L) \geq 1$, we replace $S$ by
$S/L'$, where $L'$ is the $S$--ideal generated by a maximal regular sequence in $L$.
So we may assume that $\dim R=\dim S$, and the conclusion follows by Proposition \ref{gor} and Remark \ref{e1remark}.\QED

\medskip

We give now a family of
examples based on a method of
\cite{GoNi03}.

\begin{example}{\rm Let $(S,\mathfrak{m})$ be a regular local ring of
dimension four, with an infinite residue field. Let $P_1, \ldots,
P_r$ be a family of codimension two Cohen-Macaulay ideals such that
for $i\neq j$, $P_i+P_j$ is an $\mathfrak{m}$-primary ideal.
    Define $R=S/\cap_iP_i$.

 Consider the exact
sequence of $S$-modules
\[ 0 \rar R \lar \bigoplus_i S/P_i \lar L \rar 0.\] Note that $L$ is a
module of finite support; it may be identified with
$H_{\mathfrak{m}}^1(R)$.
Let $J=(a,b)$ be an ideal of $R$  forming a system
of parameters, contained in the annihilator of $L$.\footnote{We thank
Jugal Verma for this observation.} We can assume that $a,b\in S$ form a regular
sequence in each $S/P_i$.
We are going to determine $e_1(J).$

For each integer $n$,
tensoring by $S/(a,b)^n$ we get the exact sequence
\begin{small}
\[ 0\rar \Tor_1^S(L, S/(a,b)^n) \rar R/(a,b)^n \rar \bigoplus_i S/(P_i, (a,b)^n)
 \rar L \otimes_S S/(a,b)^n \rar 0.\]
\end{small}
For $n\gg 0, $ $(a,b)^nL=0$, so we have
 that  $L \otimes_S S/(a,b)^n = L$ and
$ \Tor_1^S(L, S/(a,b)^n)= L^{n+1}$, from
 the Burch-Hilbert $(n+1)\times n$ resolution of the ideal
$(a,b)^n$. Since the  $R_i=S/P_i$  are
Cohen-Macaulay, we obtain the following Hilbert-Samuel polynomial:
\[e_0(J){{n+2}\choose {2}} -e_1(J) {{n+1}\choose{1}} + e_2(J) =
(\sum_{i=1}^r e_0(JR_i))
 {{n+2}\choose{2}}  +
(n+2)\lambda(L)-\lambda(L).\]
 It gives
\begin{eqnarray*}
e_0(JR) &=& \sum_{i=1}^r e_0(JR_i) , \\
e_1(JR) & = & -\lambda(L),\\
e_2(JR) &=& 0.
\end{eqnarray*}

}
\end{example}

\section{Embedding into balanced big Cohen-Macaulay modules}

Let $(R, \mathfrak{m})$ be a Noetherian local domain. If $R$ has a
big Cohen-Macaulay module $E$, any nonzero element of $E$ allows for an
embedding $R\hookrightarrow E$. In fact, one may assume that $E$ is a
 balanced big Cohen-Macaulay module (see \cite[Section 8.5]{BH} for a discussion).
According to the results of Hochster, if $R$
 contains a field, then there is a balanced big Cohen-Macaulay
$R$-module $E$ (\cite[8.4.2]{BH}).

To use the argument in \cite[Theorem 3.2]{chern}, in the
exact sequence
\[ 0 \rar R \lar E \lar C\rar 0,\]
we should, given any parameter ideal $J$ of $R$, pick an element
superficial for the purpose of building $\gr_J(R)$ (if $\dim R>2$)
and not contained in any associated prime of $C$ different from
$\mathfrak{m}$.
This is possible if the cardinality of the residue field is larger
than the cardinality of $\Ass (C)$.

\begin{theorem} \label{e1withBBCM} Let $(R,\mathfrak{m})$ be a
Noetherian local integral domain that is not Cohen-Macaulay and let $E$ be a balanced big
Cohen-Macaulay module. If the residue field of $R$ has cardinality
larger than the cardinality of a generating set for $E$, then
  $e_1(J)<0$
 for  any parameter ideal $J$.
\end{theorem}

Let
 $\mathbf{X}$ be a
set of indeterminates of larger cardinality than $\Ass(C)$, and
consider $R(\XX) = R[\XX]_{\mathfrak{m}R[\XX]}$. This is a Noetherian
ring (\cite{GH79}), and we are going to argue that if $E$ is a balanced big Cohen-Macaulay
$R$-module, then
  $R(\XX)\otimes_RE$ is a balanced big Cohen-Macaulay module over $R(\XX)$.
S. Zarzuela has kindly pointed out to us the following result:

\begin{theorem}[{\cite[Theorem 2.3]{Zarzuela}}]
Let $A \rightarrow B$ be a flat morphism of local rings $(A,
\mathfrak{m}),  (B, \mathfrak{n})$
and $M$ a balanced big Cohen-Macaulay $A$-module. Then, $M \otimes_A B$
is a balanced big Cohen-Macaulay $B$-module if and only if the following
two conditions hold:
\begin{itemize}

\item[{\rm (i)}]  $\mathfrak{n}(M\otimes_A B) \neq  M \otimes_A B$ and

\item[{\rm (ii)}]  For any prime ideal $\mathfrak{q} \in
\supp_B(M\otimes _A B)$,  $(1)$
$\height(\mathfrak{q} / \mathfrak{p}B) = \depth
(C_{\overline{\mathfrak{q}}})$  and $(2)$ $\height(\mathfrak{q}) +
\dim(B/\mathfrak{q}) = \dim(B)$.

\end{itemize}
Here, we denote by $\supp_A(M)$ $($small support$)$ the set of prime ideals
in $A$ with at least one non-zero Bass number in the $A$-minimal
injective resolution of $M$, $\mathfrak{p} =
\mathfrak{q}\cap A,$  $C = B/\mathfrak{p}B$ and
$\overline{\mathfrak{q}}=\mathfrak{q}C$.

Moreover, if $\mathfrak{q} \in \supp_B(M\otimes _AB)$ then
$\mathfrak{p} \in \supp_A(M)$, and
$\height(A/\mathfrak{p}) + \dim(A/\mathfrak{p}) = \dim(A)$.
\end{theorem}

\begin{corollary} Let $(R, \mathfrak{m})$ be a universally catenary
integral domain and let $E$ be a balanced big Cohen-Macaulay $R$-module. For any set $\XX$
of indeterminates and $B=R(\XX)= R[\XX]_{\mathfrak{m}[\XX]}$,
$B\otimes_RE$ is a balanced big Cohen-Macaulay $B$-module.
\end{corollary}

\smallskip

\begin{theorem} \label{e1Hochster} Let $(R, \mathfrak{m})$ be a universally catenary
integral domain containing a field. If $R$ is not Cohen-Macaulay,
then $e_1(J)<0$ for any parameter ideal $J$.
\end{theorem}

\demo Let $E$ be a balanced big Cohen--Macaulay $R$--module (\cite[8.4.2]{BH}) and consider the exact sequence
$0 \rar R \rar E \rar C \rar 0$.
Let $\mathbf{X}$ be a set of indeterminates of larger cardinality than $\Ass(C)$, and
let $R(\XX) = R[\XX]_{\mathfrak{m}R[\XX]}$. By applying Theorem~\ref{e1withBBCM} to
the exact sequence $0 \rar R(\XX) \rar E \otimes R(\XX) \rar C \otimes R(\XX) \rar 0$, the assertion is proved.
\QED

\section{Filtered modules}

The same relationship discussed above between the signature of $e_1(J)$
and the Cohen-Macaulayness of $R$ holds true when modules are examined.
Recall that if a Noetherian local ring $R$ is embedded into either a maximal Cohen--Macaulay module (\cite[Theorem 3.1]{chern}) or a balanced big Cohen--Macaulay module (\cite[proof of Theorem 3.2]{chern}), then whenever $R$ is not Cohen--Macaulay, we have $e_1(J) <0$ for any parameter ideal $J$. Now we use the same arguments as in \cite[Theorems 3.1, 3.2]{chern} in order to extend the validity of Theorem~\ref{e1dMCM}
in the following manner.

\begin{theorem} \label{e1formodules}
Let $(R, \mathfrak{m})$ be a Noetherian local ring of dimension
$d\geq 1$
and let $M$ be a finitely generated module embedded in a
maximal Cohen--Macaulay module $E$.  Then $M$ is Cohen-Macaulay if and only if
$e_1(J, M) \geq 0$ for any ideal $J$ generated by a system of
parameters of $M$.
\end{theorem}

A variation that uses Theorem~\ref{e1Hochster} is the following.

\begin{theorem} \label{e1Hochster4modules} Let $(R, \mathfrak{m})$ be a universally catenary
integral domain containing a field and let $M$ be a finitely generated torsion--free
$R$-module. If $M$ is not Cohen-Macaulay,
then $ e_1(J, M)<0$ for any ideal $J$ generated by a system of
parameters of $M$.
\end{theorem}

\demo By assumption $M$ is a submodule of a finitely generated free
$R$-module, which can be embedded into a
finite direct sum of balanced big Cohen-Macaulay modules.
The argument of \cite[Theorem 3.2]{chern} applies again.\QED

\bigskip

Let now $R=k[x_1, \ldots, x_d]$ be a ring of polynomials over the field
$k$, and let $M$ be a finitely generated graded $R$-module. Suppose
$\dim M=d$. For $J=(x_1, \ldots, x_d)$ we can apply
Theorem~\ref{e1formodules} to $M$ in a manner that uses the
Hilbert--Samuel
function information of the native grading of $M$.

\begin{theorem} \label{e1vsCM0}
Let $R=k[\xx]=k[x_1, \ldots, x_d]$, $d\geq 2$,  be a ring of polynomials over the field
$k$, and let $M$ be a finitely generated graded $R$-module generated
in degree $0$. If $M$ is torsion--free, then $M$
 is a free $R$-module if and only if $e_1(M)=0$.
\end{theorem}

\demo Since $M$ is generated in degree $0$, $M \simeq \gr_{\xx}(M)$.
By assumption, $M$ can be embedded in a free $R$-module $E$ (not
necessarily by a homogeneous homomorphism).  Now we apply
Theorem~\ref{e1formodules}. \QED

\bigskip

If $M$ is generated in degree $a>0$, we have the equality
\[ \lambda(M/(\xx)^{n+1}M) = \sum_{k=0}^n \lambda(M_{a+k}),\]
so the Hilbert coefficients satisfy
\begin{eqnarray*}
e_0(\gr_{\xx}(M)) &=& e_0(M[a])= e_0(M),\\
e_1(\gr_{\xx}(M)) &=& e_1(M[a])= e_1(M)-ae_0(M).
\end{eqnarray*}

\begin{corollary} \label{e1vsCMa}
Let $R=k[\xx]=k[x_1, \ldots, x_d]$, $d\geq 2$,  be a ring of polynomials over the field
$k$, and let $M$ be a finitely generated graded $R$-module generated
in degree $a\geq 0$. If $M$ is torsion--free, then $e_1(M)\leq
ae_0(M)$, with equality if and only if
$M$
 is a free $R$-module.
\end{corollary}

\newpage

\end{document}